\documentclass[12pt]{amsart}  

\newenvironment{pf}{\emph{Proof. }}{\qed\smallskip}
\newtheorem{thm}[subsection]{Theorem}
\newtheorem{cor}[subsection]{Corollary}
\newtheorem{prop}[subsection]{Proposition}
\newtheorem{lemma}[subsection]{Lemma}

\theoremstyle{remark}
\newtheorem{remark}[subsection]{Remark}

\theoremstyle{definition}
\newtheorem{defn}[subsection]{Definition}
\newtheorem{que}[subsection]{Question}

\numberwithin{equation}{section}

\newcommand{\bbA}{\mathbb{A}}
\newcommand{\bbC}{\mathbb{C}}

\newcommand{\bbN}{\mathbb{N}}
\newcommand{\bbZ}{\mathbb{Z}}
\newcommand{\bbP}{\mathbb{P}}
\newcommand{\bbQ}{\mathbb{Q}}

\newcommand{\cM}{\mathcal{M}}
\newcommand{\cV}{\mathcal{V}}
\newcommand{\cH}{\mathcal{H}}

\newcommand{\sign}{\operatorname{sign}}

\newcommand{\Ker}{\operatorname{Ker}}

\renewcommand{\dim}{\text{dim}}

\newcommand{\printname}[1]
  {\smash{\makebox[0pt]{\hspace{-2.0in}\raisebox{8pt}{\tiny #1}}}}

\thanks
{The first author was supported in part by the NSF grant DMS-0100537. 
The second 
author was supported in part by the NSA grant MDA904-01-1-0020 and the 
CRDF grant RM1-2089}

\begin{document}

\title{Motivic measures and stable birational geometry}

\author{Michael Larsen}
\address{Department of Mathematics, Indiana University,
Bloomington, IN 47405, USA}
\email{larsen@math.indiana.edu}
\author{Valery A.~Lunts}
\address{Department of Mathematics, Indiana University,
Bloomington, IN 47405, USA}
\email{vlunts@indiana.edu}

\begin{abstract} We study the motivic Grothendieck group of algebraic 
varieties from the point of view of stable birational geometry. In 
particular, we obtain a counter-example to a conjecture of M.~Kapranov 
on the rationality of motivic zeta-function.
\end{abstract}

\maketitle

\section{Introduction}

\subsection{Grothendieck ring of varieties}
Fix a field $k$. Let $\cV _k$ denote the category of $k$-varieties. 
By a $k$-variety we mean a reduced separated scheme of finite type over $k$. 
Consider the Grothendieck ring $K_0[\cV _k]$: this is the abelian group 
generated by isomorphism classes of $k$-varieties, subject to the relations 
$[X-Y]=[X]-[Y]$, where $Y$ is a closed in $X$. The product over $k$ turns 
it into a commutative ring with 1. It appears that very little is known about 
this interesting ring. For example, one would like to know the answer to the 
following basic question.

\medskip

\begin{que} Let $X$,$Y$ be $k$-varieties such that $[X]=[Y]$. Is it 
possible to partition $X$ and $Y$ by a finite number of locally closed 
subvarieties which are pairwise isomorphic?
\end{que}

\subsection{Motivic measures} Let $A$ be a commutative ring. An $A$-valued 
{\it motivic measure} on $\cV _k$ is a ring homomorphism 
$$\mu :K_0[\cV _k]\rightarrow A.$$
The case $A=K_0[\cV _k]$, $\mu =id$ is the universal motivic measure.

\subsection{Motivic zeta-function} For a $k$-variety $X$ denote by $X^{(n)}$ 
the $n$-fold symmetric product of $X$. Given a motivic measure 
$\mu :K_0[\cV _k]\rightarrow A$ we define (following [Ka]) the 
{\it motivic zeta-function} of $X$ as the formal power series
$$\zeta _{\mu }(X,t)=1+\sum_{n=1}^\infty\mu([X^{(n)}])t^n\in A[[t]].$$

\medskip

For example,
$$\zeta_{id}(\bbP^1,t)=\sum_{n=0}^\infty[\bbP^n]t^n=\frac{1}
{(1-t)(1-[\bbA^1]t)},$$

$$\zeta_{id}(E,t)=1+\sum_{n=1}^\infty [E][\bbP^{n-1}]t^n=
\frac{1+([E]-[\bbP^1])t+[\bbA^1]t^2}{(1-t)(1-[\bbA^1]t)}$$
for any elliptic curve $E$, and

$$\zeta_{id}(\bbP^2,t)=\frac{1}{(1-t)(1-[\bbA^1]t)(1-[\bbA^2]t)}.$$

Checking these formulas requires a certain amount of cutting and pasting, using the
relation $[X] = [X-Y]+[Y]$ systematically.

\medskip

Kapranov proves the following theorem in [Ka]:

\begin{thm} Let $A$ be a field and $\mu :K_0[\cV _k]\rightarrow A$ 
be a motivic measure. Let $X \in \cV _k$ be a curve. Then 
$\zeta _{\mu}(X,t)$ is a rational function.
\end{thm}

In the same paper Kapranov remarks that it is natural to expect a similar 
rationality result for varieties of higher dimension. We give a negative 
solution to this problem.

\begin{thm} Assume that $k=\bbC$. There exists a field $\cH$ and a 
motivic measure $\mu :K_0[\cV _\bbC]\rightarrow \cH$ with the following 
property: if $X$ is a smooth complex projective surface  
such that $P_g(X)=h^{2,0}(X)\geq 2$, then the zeta-function
$\zeta _{\mu}(X,t)$ is not rational. 
\end{thm}

We expect that a similar result holds for any smooth projective variety of 
even dimension and positive Kodaira dimension.

\subsection{Stable birational geometry}
The above theorem on the irrationality of the zeta-function follows easily 
from our analysis of the Grothendieck ring $K_0[\cV _\bbC]$ from the point 
of view of stable birational equivalence of varieties. Namely, recall that 
(irreducible) varieties $X$,$Y$ are stably birational if $X\times \bbP ^k$ 
is birational to $Y\times \bbP ^l$ for some $k,l\geq 0$. Let $SB$ denote the 
multiplicative monoid of classes of stable birational equivalence of 
varieties; let $\bbZ [SB]$ be the corresponding monoid ring. Consider the 
principal ideal $I\subset K_0[\cV _\bbC]$ generated by the class of 
the affine line $\bbA ^1$. The main result of the next section 
 is the following 
isomorphism of rings.  
$$K_0[\cV _\bbC]/I\simeq \bbZ [SB].$$
This isomorphism is interesting in its own right. It implies, in particular, 
that classes of stable birational equivalence from a $\bbZ$-basis of the 
group $K_0[\cV _\bbC]/I$. The hard part in establishing this 
isomorphism is to show that the map from the LHS to the RHS is well defined. 
This is essentially the content of Theorem \ref{mt}.

We also show that for a ring homomorphism $\phi :K_0[\cV _\bbC]\to B$ the 
following conditions are equivalent

i) $I\subset \Ker (\phi )$,

ii) if $X$, $Y$ are smooth complete varieties which are birational, then 
$\phi ([X])=\phi ([Y])$.

\section{Grothendieck ring of varieties and stable birational equivalence}

\subsection{The Grothendieck ring of varieties revisited} We make a few 
remarks about the Grothendieck ring $K_0[\cV _k]$ which will be used later. 

\medskip

\noindent A. If a variety $X$ is partitioned by locally closed subvarieties 
$X_1,...X_n$, then $$[X]=\sum [X_i].$$

\medskip

\noindent B. Every variety can be partitioned by a finite number of 
smooth varieties. Hence classes of smooth varieties generate the group 
$K_0[\cV_k]$. 

\medskip

\noindent C. Let $E\to X$ be a vector bundle of dimension $d$. 
Then it is locally trivial 
in the Zariski topology. Hence the same is true about its projectivization
$\overline{E}\to X$. Thus $[\overline{E}]=[X][\bbP ^{d-1}]$ in the ring 
$K_0[\cV _k]$.

\medskip

\noindent D. Let $f:X\to Y$ be a proper morphism of smooth varieties, 
which is a blowup with a smooth center $Z\subset Y$ of codimension $d$. 
Then the projective bundle $f^{-1}(Z)\to Z$ is the projectivization of the 
normal bundle of $Z$ in $Y$. Hence, by the previous remark
 $[f^{-1}(Z)]=[Z][\bbP^{d-1}]$ in the ring 
$K_0[\cV _k]$. 

\medskip

\noindent E. Assume that $char(k)=0$. By Hironaka's theorem every smooth 
variety $X$ is isomorphic to a dense open subset of a smooth complete 
variety. Hence the group $K_0[\cV_k]$ is generated by smooth complete 
varieties.

\medskip

Let $k=\bbC$. We will be using the following theorem of Wlodarczyk (see [W],
[AKMW]). 

\begin{thm} \label{Wlod}
Let $\phi :X_1\to X_2$ be a rational birational map of 
smooth complete varieties. Let $U\subset X_1$ be an open subset where 
$\phi $ is an isomorphism. Then 
$\phi =\phi _l{\scriptstyle\circ}\cdots{\scriptstyle\circ} \phi _1$, where each $\phi _i$ is a rational 
birational map of smooth complete varieties 
and either $\phi _i$ or 
$\phi ^{-1}_i$ 
is a regular map which is a blowup with a smooth center disjoint 
from $U$. 
\end{thm}

The basic result of this paper is the following:

\begin {thm}
\label{mt}
Put $k=\bbC$. Let $G$ be an abelian commutative monoid and $\bbZ [G]$ 
be the corresponding monoid ring. Denote by $\cM$ the multiplicative monoid 
of isomorphism classes of smooth complete irreducible varieties. Let 
$$\Psi :\cM \rightarrow G$$
be a homomorphism of monoids such that 

(i) $\Psi ([X])=\Psi ([Y])$ if $X$ and $Y$ are birational;

(ii) $\Psi ([\bbP ^n])=1$ for all $n\geq 0$. 

Then there exists a unique ring homomorphism 
$$\Phi :K_0[\cV _\bbC]\rightarrow \bbZ[G]$$
such that $\Phi ([X])=\Psi ([X])$ for $[X]\in \cM$.
\end{thm}

\begin{pf} To simplify notation we will write $\Psi (X)$ and $\Phi (X)$ for 
$\Psi ([X])$ and $\Phi ([X])$ respectively. 

We will define the elements $\Phi (X)\in \bbZ [G]$ by induction on the 
dimension of the variety $X$. The induction step will require checking that 

i) $\Phi (X)$ is well defined for each variety $X$;

ii) $\Phi $ preserves the defining relations of the abelian group 
$K_0[\cV _\bbC]$;

iii) $\Phi $ is multiplicative.

Let us formulate a series of constructions and assertions which depend on $n$. 

\noindent{\bf Construction} $A_n$. If $X$ is an irreducible smooth complete variety of 
dimension 
$\leq n$ then $\Phi (X):=\Psi (X)$.

\noindent{\bf Construction} $B_n$. 
Let $X$ be a smooth variety of dimension $\leq n$ 
with connected components $X_1,...,X_k$. For each $i$ choose an open 
embedding $X_i\hookrightarrow \overline{X_i}$ where $\overline{X_i}$ is 
smooth, complete, and irreducible. Then put 
$$\Phi (X):=\sum \Phi (\overline{X_i})-\sum \Phi (\overline{X_i}-X_i).$$

\noindent{\bf Construction} $C_n$. Let $X$ be an arbitrary variety of dimension $\leq n$. 
Then put $$\Phi (X):=\Phi (X-X^{sing})+\Phi (X^{sing}).$$

\medskip

\noindent{\bf Assertion} $D_n$. Let $X$, $Y$ be varieties of dimension 
$\leq n$, and $f:X\to Y$ be a morphism with the following property:
 there exists a stratification of $Y$ by locally closed subvarieties 
$Y_i$, so that $f^{-1}(Y_i)\simeq \bbP ^{n_i}\times Y_i$ for some $n_i\geq 0$ 
and $f:f^{-1}(Y_i)\to Y_i$ is the projection. Then 
$$\Phi (X)=\Phi (Y).$$ 

\medskip

\noindent{\bf Assertion} $E_n$. Let $X$ be a variety of dimension $\leq n$ and 
$Y\subset X$ a closed subvariety. Then 
$$\Phi (X)=\Phi (Y) +\Phi (X-Y).$$

\medskip

\noindent{\bf Assertion} $F_n$. Let $X$,$Y$ be varieties such that 
$\dim X +\dim Y\leq n$, then 
$$\Phi (X\times Y)=\Phi (X)\cdot \Phi (Y).$$

\medskip

We will use the following logic in proving the theorem. Assume that $\Phi $ 
is constructed according to Constructions $A_{n-1}$,  $B_{n-1}$, $C_{n-1}$,  
Construction $B_{n-1}$ is unambiguous, and Assertions $D_{n-1}$, 
$E_{n-1}$, $F_{n-1}$ are 
proved. Then, in particular, $\Phi $ is defined on all classes $[X]$ of 
varieties $X$ of dimension $\leq n-1$ in such a way that 

1) $\Phi (X)=\Psi (X)$ if $X$ is irreducible, smooth, and complete.

2) $\Phi (X)=\Phi (Y)+\Phi (X-Y)$ if $Y\subset X$ is a closed. 

3) $\Phi (X\times Y)=\Phi (X)\cdot \Phi (Y)$ if $\dim X+\dim Y\leq n-1$. 

We extend $\Phi$ by linearity to linear combinations of such classes. Then 
we use Definitions $A_n$, $B_n$, $C_n$ to extend 
$\Phi $ to classes of varieties 
of dimension $\leq n$ and prove Assertions $D_n$, $E_n$, $F_n$ 
which are needed to 
ensure that Definition $B_n$ is unambiguous and that this extension of $\Phi$ 
satisfies the properties 1),2),3) above. This will prove the existence 
statement of the theorem. The uniqueness is clear since classes of smooth 
complete irreducible varieties generate the group $K_0[\cV _\bbC]$. 

\medskip

\noindent{\it Base case.} We can use Definitions $A_0$, 
$B_0$, $C_0$ without ambiguity, and Assertions $D_0$, $E_0$, $F_0$ 
obviously hold.

\medskip

\noindent{\it Induction step $(n-1)\to n$.} 
Assume that the map $\Phi$ has been defined using  Constructions $A_{n-1}$, 
 $B_{n-1}$, $C_{n-1}$, Construction $B_{n-1}$ is unambiguous, 
and Assertions $D_{n-1}$, $E_{n-1}$, $F_{n-1}$ are true. 

Let $X$ be an irreducible smooth complete variety of dimension $n$. 
Then define $\Phi (X)$ according to $A_n$. 

Let $X$ be a smooth $n$-dimensional variety. Use $B_n$ to define 
$\Phi (X)$.  We prove that it is independent of the choice of smooth 
compactifications (such compactifications exist by Hironaka's theorem). 
We may assume that $X$ is irreducible. Let $X\hookrightarrow \overline{X}$, 
$X\hookrightarrow \overline{X}^\prime$ be two open embeddings
with $\overline{X}$, $\overline{X}^\prime$  smooth and complete. 
Put $Y:=\overline{X}-X$, $Y^\prime:=\overline{X}^\prime$. We need to show that 
$\Phi (\overline{X})-\Phi (Y)=\Phi (\overline{X}^\prime)-\Phi (Y^\prime)$. 
Since $\Psi $ is a birational invariant, Construction $A_n$ implies 
$\Phi (\overline{X})=\Phi (\overline{X}^\prime)$. So it remains to show that 
$\Phi (Y)=\Phi (Y^\prime)$. By the Theorem~\ref{Wlod} we may reduce
to the case that 
there exists a morphism $f:\overline{X}^\prime \to \overline{X}$ which is a 
blowup with a smooth center $Z\subset Y$. Then by Remark C above the 
map $f:f^{-1}(Z)\to Z$ is a Zariski locally trivial fibration with fibre 
$\bbP^k$. So by $D_{n-1}$, $\Phi (Y)=\Phi (Y^\prime)$. This justifies the 
definition $B_n$. 

Now for a general $n$-dimensional variety $X$ we define $\Phi (X)$ according 
to $C_n$. 

It remains to prove the Assertions $D_n$, $E_n$, $F_n$. 
For a variety $W$ put $W^{ns}:=W-W^{sing}$. 

\medskip

\noindent{\it Proof of} $E_n$. 
Let $X$ be a variety of dimension $n$ and $Y\subset X$ be a closed subvariety. 
Put $U:=X-Y$. We need to prove that $\Phi (X)=\Phi (U)+\Phi (Y)$. 

\medskip

Assume first that $X$ is {\it smooth}. If $W$ is any variety   and 
$W_1,...,W_s$ are the connected components of $W$, then by Constructions $B_n$ and
$C_n$, $\Phi (W)=\sum \Phi (W_i)$. So we may assume that X is {\it connected}.
If $\dim Y=n$, then $U=\emptyset$ and we are done. So assume that $\dim Y<n$. 
Let $X\hookrightarrow \overline{X}$ be a smooth compactification of $X$, 
$Z:=\overline{X}-X$, $S:=\overline{X}-U=Z\coprod Y$. Then by Construction $B_n$,
$$\Phi (X)=\Phi (\overline{X})-\Phi (Z),\quad \Phi (U)=\Phi (\overline{X})-
\Phi (S).$$
By $E_{n-1}$ $\Phi(S)=\Phi (Z)+\Phi (Y)$. So 
$$\Phi (X)=\Phi (U)+\Phi(Y).$$

\medskip

Now let $X$ be any variety. We have 
$$U^{sing}=X^{sing}\cap U,\quad U^{ns}=X^{ns}\cap U.$$
By $C_n$, 
$$\Phi(X)=\Phi (X^{ns})+\Phi (X^{sing}),$$
$$\Phi(U)=\Phi (U^{ns})+\Phi (U^{sing}).$$
By the argument above,
$$\Phi (X^{ns})=\Phi (U^{ns})+\Phi (Y\cap X^{ns}),$$
and therefore
$$\Phi (X)=\Phi (U)-\Phi (U^{sing})+\Phi (Y\cap X^{ns})+\Phi (X^{sing}).$$
So it remains to prove that 
$$\Phi (Y)=\Phi (X^{sing})+\Phi (Y\cap X^{ns})-\Phi (U^{sing}).$$
Since $X^{sing}=U^{sing}\coprod (Y\cap X^{sing})$, by $E_{n-1}$, 
$$\Phi (X^{sing})=\Phi (U^{sing})+\Phi (Y\cap X^{sing}).$$
So it suffices to 
prove that $\Phi(Y)=\Phi (Y\cap X^{ns})+\Phi (Y\cap X^{sing})$. The last 
equality is $E_n$ with $X,Y,U$ replaced by $Y,Y\cap X^{sing},Y\cap X^{ns}$. 
Since we may assume that $Y$ is a proper subset of $X$ the proof is finished 
by Noetherian induction on $X$. 

\medskip

\noindent{\it Proof of} $F_n$. Let $X,Y$ be varieties such that 
$\dim X+\dim Y=n$. We need to prove that $\Phi (X\times Y)=\Phi (X)\cdot 
\Phi (Y)$. If $\dim Y=0$, then we are done. So may assume that 
$\dim X,\dim Y\geq1$. By $E_n$ we have 
\begin{eqnarray*}
\Phi(X\times Y)&=&\Phi (X^{ns}\times Y^{ns})+\Phi (X^{ns}\times Y^{sing})\\
&&\qquad+\Phi (X^{sing}\times Y^{ns})+\Phi (X^{sing}\times Y^{sing}),\\
\end{eqnarray*}

\begin{eqnarray*}
\Phi(X)\cdot \Phi(Y)&=&\Phi (X^{ns})\cdot \Phi (Y^{ns})+\Phi (X^{ns})\cdot 
\Phi(Y^{sing})\\
&&\qquad+
\Phi (X^{sing})\cdot \Phi(Y^{ns})+\Phi (X^{sing})\cdot \Phi(Y^{sing}).\\
\end{eqnarray*}

Hence by $F_{n-1}$ we may assume that $X,Y$ are smooth (and connected). 
Let $X\hookrightarrow \overline{X}$, $Y\hookrightarrow \overline{Y}$ be 
smooth compactifications of $X,Y$, $X^\prime :=\overline{X}-X$, 
$Y^\prime :=\overline{Y}-Y$. Again by $E_n$,
$$\Phi(\overline{X})\cdot \Phi( \overline{Y})=
\Phi (X)\cdot \Phi (Y)+\Phi (X)\cdot 
\Phi(Y^\prime)+
\Phi (X^\prime)\cdot\Phi(Y)+\Phi (X^\prime)\cdot \Phi(Y^\prime),$$
$$\Phi(\overline{X}\times \overline{Y})=\Phi (X\times Y)+
\Phi (X\times Y^\prime)+
\Phi (X^\prime \times Y)+\Phi (X^\prime\times Y^\prime).$$
We have $\Phi (\overline{X}\times \overline{Y})=
\Phi(\overline{X})\cdot \Phi(\overline{Y})$, since 
$\Psi$ is multiplicative. Now using $F_{n-1}$ we conclude that 
$\Phi (X\times Y)=\Phi (X)\cdot \Phi (Y).$

\medskip

\noindent{\it Proof of} $D_n$. Using $E_n$ repeatedly we may assume that 
$X=Y\times \bbP^k$ and $f:X\to Y$ is the projection. Then by $F_n$,
$\Phi (X)=\Phi (Y)\cdot \Phi (\bbP^k)$, and $\Phi (\bbP ^k)=
\Psi (\bbP ^k)=1$. 
\end{pf}

\begin{remark} The Grothendieck group $K_0[\cV_{\bbC}]$ is generated by 
classes of smooth complete varieties. In [L] Looijenga asserts
that it suffices to consider the following relations: let $X,Y$ be smooth 
complete and $f:X\to Y$ be a morphism which is a blowup with a smooth center 
$Z\subset Y$; then 
$$[X]-[f^{-1}(Z)]=[Y]-[Z].$$
This is a strong result, which immediately implies our Theorem \ref{mt}. 
Indeed, $f^{-1}(Z)$ is birational to $Z\times \bbP^k$. So by the hypotheses 
of  Theorem~\ref{mt}, $\Psi (X)=\Psi (Y)$, $\Psi (f^{-1}(Z))=
\Psi (Z\times \bbP ^k)=\Psi (Z)$. Therefore the desired ring homomorphism 
$\Phi $ exists. However, since we could not produce a proof of Looijenga's 
result, we chose to give an argument that does not depend on it. 
(Looijenga informs us that 
a proof is being written up by one of his students.) 
\end{remark}

\subsection{The universal homomorphism $\Psi :\cM \to G$}
There exists a universal homomorphism of monoids $\Psi$ which satisfies 
the hypotheses of Theorem \ref{mt}. Namely, recall that varieties $X,Y$ are 
{\it stably birational} if $X\times \bbP^k$ and $Y\times \bbP^l$ are 
birational for some $k,l\geq 0$. Denote by $SB$ the multiplicative monoid 
of stable birational equivalence classes of varieties. We have a tautological 
(surjective) homomorphism $\Psi _{SB}:\cM \to SB$ which satisfies the 
hypotheses (i),(ii) of Theorem \ref{mt} (with $\Psi =\Psi _{SB}$, $G=SB$). By 
definition any homomorphism $\Psi $ as in the theorem factors through 
$\Psi _{SB}$. 
Denote by $\Phi _{SB}:K_0[\cV_{\bbC}]\to \bbZ [SB]$ 
the ring homomophism corresponding to $\Psi _{SB}$ by the theorem.
We obtain the following immediate corollary.

\begin{cor} Let $X_1,...X_k$, $Y_1,...Y_m$ be smooth complete varieties.
 Let $m_i,n_j\in \bbZ$ be such that 
$$\sum m_i[X_i]=\sum n_j[Y_j]$$
in $K_0[\cV _{\bbC}]$. Then $k=m$ and after renumbering the varieties 
$X_i$ and $Y_i$ are stably birational and $m_i=n_i$. 
\end{cor}

\begin{pf} 
Applying the ring homomorphism $\Phi _{SB}$ to the above equality we obtain 
the equality in the monoid ring $\bbZ[SB]$:
$$\sum m_i\Psi_{SB}(X_i)=\sum n_j \Psi _{SB}(Y_j)$$
and the proposition follows.
\end{pf}

The above corollary means that any variety is a {\it unique} (up to a stable 
birational equivalence) linear combination (in $K_0[\cV _{\bbC}]$) of smooth 
complete varieties. This is in the spirit of the basic Question formulated 
in the introduction. The difference is, of course, that instead of cutting 
varieties in pieces we complete and resolve singularities.

The next proposition clarifies the relation between the 
Grothendieck ring $K_0[\cV _{\bbC}]$ and the monoid ring $\bbZ[SB]$. 

\begin{prop} The kernel of the (surjective) homomorphism 
$\Phi _{SB}:K_0[\cV_{\bbC}]\to \bbZ[SB]$ is the principal ideal generated 
by the class $[\bbA ^1]$ of the affine line $\bbA ^1$. 
\end{prop}

\begin{pf} Since $\Phi _{SB}([\bbP ^1])=\Phi _{SB}([1]+[\bbA ^1])=1$, we have 
$\Phi _{SB}([\bbA ^1])=0$.

Let $a\in \Ker (\Phi _{SB})$. Express $a$ as a linear combination
$$a=[X_1]+...+[X_k]-[Y_1]-...-[Y_l],$$
where $X_i$, $Y_j$ are smooth and complete. Since 
$$\Phi _{SB}(a)=\sum \Psi _{SB}(X_i)-\sum \Psi _{SB}(Y_j)=0,$$
we get $k=l$ and, after renumbering, $X_i$ is stably birational to $Y_i$. 
Thus it suffices to show that if $X$, $Y$ are smooth, complete, and stably 
birational, then $[X]-[Y]\in K_0[\cV _{\bbC}]\cdot[\bbA ^1]$. Note that 
$$[X\times \bbP^k]-[X]=[X]\cdot [\bbA ^1+\bbA^2+...+\bbA ^{k}],$$
so we may 
assume that $[X]$ and $[Y]$ are birational. Moreover by Theorem~\ref{Wlod} we may 
assume that $X$ is a blowup of $Y$ with a smooth center $Z\subset Y$ and 
exceptional divisor $E\subset X$. Then $[E]=[\bbP ^t]\cdot [Z]$ for some 
$t$ and 
$$[X]-[Y]=[E]-[Z]=([\bbA ^1]+[\bbA^2]+\cdots+[\bbA ^{t}])\cdot [Z].$$
\end{pf}

The next proposition ``explains'' the role of $\bbA ^1$ in birational 
geometry.

\begin{prop}
Let $\alpha :K_0[\cV _\bbC]\to B$ be a ring homomorphism (i.e. 
$\alpha $ is a motivic measure). The following conditions are equivalent:

i) $\alpha ([\bbA ^1])=0$,

ii) if smooth complete varieties $X$, $Y$ are birational, then 
$\alpha ([X])=\alpha ([Y])$.

If these conditions hold, then $\alpha ([Z])=\alpha ([W])$ for any 
smooth complete varieties $Z$, $W$ which are stably birationally 
equivalent
\end{prop}

\begin{pf} Assume i). Then $\alpha $ factors through the homomorphism 
$\Phi _{SB}$ and ii) follows. Assume ii). Let $\tilde{X}\to X$ be a blowup 
of a smooth complete surface at a point. Then $[\tilde{X}]=[X]+[\bbA ^1]$ and 
hence $\alpha ([\bbA ^1])=0$. 

To prove the last assertion note that 
$\alpha ([\bbA^1])=0$ implies that $\alpha ([\bbP ^n])=1$. 
\end{pf}

The last assertion of the proposition means that birational motivic 
measures are automatically stably birational.

\section{Irrationality of the zeta-function}

\subsection{The motivic measure $\mu _h:K_0[\cV _{\bbC}]\to \cH$} 
Let $C\subset \bbZ [t]$ be the multiplicative monoid of polymonials with a
 positive leading coefficient. Consider the corresponding monoid ring 
$\bbZ [C]$.

\begin{lemma} The ring $\bbZ[C]$ is an integral domain.
\end{lemma}

\begin{pf} The ring $\bbZ[t]$ is factorial and any element of $C$ is a 
unique product of elements of $C$, which are prime in $\bbZ [t]$ (the 
only unit in $C$ is 1). Thus $C$ is isomorphic to the group $\oplus \bbN$, 
where the summation is over all prime elements of $\bbZ[t]$. Hence 
$\bbZ [C]$ is a polynomial ring, so it is an integral domain.
\end{pf}

\begin{defn} Let $\cH$ be the field of fractions of $\bbZ [C]$. 
\end{defn}

\begin{defn}
For a smooth projective irreducible complex variety $X$ 
 of dimension $d$ define 
$$\Psi_h(X):=1+h^{1,0}(X)t+...+h^{d,0}(X)t^d\in C.$$
For any smooth complete complex irreducible variety $Z$ put $\Psi _h(Z)=\Psi _h(X)$, 
where $X$ 
is a smooth projective variety which is birational to $Z$. 
It is well known that if smooth projective varieties $X$ and $Y$ 
are birational then $\Psi_h(X)=\Psi_h(Y)$ ([Ha], Ch. 2, Exercise 8.8). 
Therefore $\Psi _h$ is well defined. 
The K\"unneth formula implies that $\Psi _h$ is a homomorphism from 
the multiplicative monoid $\cM$ of isomorphism classes of smooth complete 
 irreducible 
varieties to $C$. It 
 satisfies the hypotheses of Theorem \ref{mt}. Thus it extends to 
a motivic measure 
$$\mu _h:K_0[\cV _\bbC]\rightarrow \cH.$$
\end{defn}

\begin{defn} For a smooth projective irreducible variety $X$ of dimension $d$ 
denote as usual 
$$P_g(X)=h^{d,0}(X)=h^0(X,\omega _X),$$
where $\omega _X$ is the canonical line bundle on $X$. For an arbitrary 
irreducible variety $Z$ of dimension $d$ put 
$$P_g(Z):=P_g(X),$$
where $X$ is any smooth projective variety in the birational 
class of $Z$. Thus $P_g$ becomes a multiplicative function from the 
collection of isomorphism classes of irreducible varieties to natural numbers. 
\end{defn}

\begin{lemma}\label{nocancel}
Let $Y_1,...Y_s,Z$ be irreducible varieties of dimension $d$. 
Assume that $\mu _h([Z])= \sum n_i\mu _h([Y_i])$ for some 
$n_i\in \bbZ $. If $P_g(Z)\neq 0$, then  $P_g(Z)=P_g(Y_i)$ for 
some $i$.
\end{lemma}

\begin{pf}
Note that for any irreducible variety $W$ of dimension $d$ we have the 
equality in $K_0[\cV _{\bbC}]$
$$[W]=[\overline{W}]+\sum m_j[W_j],$$
where $\overline{W}$ is a smooth projective variety in the birational 
class of $W$ and $W_j$'s are smooth projective varieties of dimension 
$< d$. Thus replacing $Y_1,...,Y_s,Z$ by smooth projective 
varieties $\overline{Y}_1,...,\overline{Y}_s,\overline{Z}$ from the same 
birational class we obtain an equality 
$$\mu _h([\overline{Z}])=\sum n_i\mu _h([\overline{Y}_i])+\sum
l_{\beta}\mu _h([X_{\beta }])$$
for some $l_{\beta}\in \bbZ$ and some smooth projective irreducible varieties 
$X_{\beta }$ of dimension $< d$. By definition this means
$$\Psi _h(\overline{Z})=\sum n_i\Psi _h(\overline{Y}_i)+\sum l_{\beta }
\Psi _h(X_{\beta}),$$
and the lemma follows. 
\end{pf}

\begin{thm} 
Consider the motivic measure $\mu _h:K_0[\cV _{\bbC}]\to \cH$. 
Let $X$ be a smooth complex projective surface. 
Assume that $P_g(X)\geq 2$. Then the zeta function 
$\zeta _{\mu _h}(X,t)$ is not rational.
\end{thm}

\begin{pf} For a variety $Y$ we will write for short $\mu (Y)=\mu _h([Y])$. 
Let $X$ be as in the theorem and assume that the zeta-function 
$$\zeta _{\mu _h}(X,t)=1+\sum_{n=1}^\infty\mu (X^{(n)})\in \cH[[t]]$$ 
is rational. Recall the following characterization of rational power series 
(with coefficients in a field). 
A power series $\sum a_it^i$ is a rational function if and only if there exist 
$n> 0$, $n_0> 0$ such that for each $m>n_0$ the determinant
of the matrix 
$$\left(
\begin{array}{cccc}
a_m & a_{m+1} & \ldots & a_{m+n} \\
a_{m+1} & a_{m+2} & \ldots & a_{m+n+1} \\
\vdots & \vdots & \ddots & \vdots \\
a_{m+n} & a_{m+n+1} & \ldots & a_{m+2n}
\end{array}
\right)
$$
is zero. 

In case $\sum a_it^i=\zeta _{\mu _h}(X,t)$ this equality implies 

\renewcommand\theequation{$*$}
\begin{equation}
\sum_{\sigma\in S_{n+1}}
\sign(\sigma) \mu \bigl(X^{(m-1+\sigma(1))}\times X^{(m+\sigma(2))}\times\cdots\times
X^{(m+n-1+\sigma(n+1))}\bigr)=0.
\end{equation}

Note that the summand $\mu \bigl(X^{(m)}\times X^{(m+2)}\times\cdots\times 
X^{(m+2n)}\bigr)$ appears exactly once.

\medskip

\noindent{\it Claim.} {\it There exists $m>0$ such that 
$$P_g\bigl(X^{(m)}\times\cdots\times 
X^{(m+2n)}\bigr)=
P_g\bigl(X^{(m-1+\sigma(1))}\times\cdots\times
X^{(m+n-1+\sigma(n+1))}\bigr)$$
implies that $\sigma$ is the identity permutation. }

\medskip

Assuming the claim we apply Lemma~\ref{nocancel} with 
$Z=X^{(m)}\times X^{(m+2)}\times\cdots\times X^{(m+2n)}$ and $Y_j$'s being 
the other summands in ($*$). We get a contradiction and the theorem follows. 
So it remains to prove the claim.

\begin{lemma} Let $r=P_g(X)$. Then 
$$P_g(X^{(n)})=\binom{r+n-1}{r-1}.$$
\end{lemma}

Assume the lemma and let us prove the claim. By our assumption 
$r= P_g(X)\geq 2$, so

\begin{eqnarray*}
\lefteqn{P_g\bigl(X^{(m-1+\sigma(1))}\times X^{(m+\sigma(2))}\times\cdots\times 
	X^{(m+n-1+\sigma(n+1))}\bigr)} \\
& & =  \prod_{j=0}^n P_g\bigl(X^{(m+j-1+\sigma(j+1))}\bigr) \\
& & =  \frac{1}{(r\!-\!1)!^{n+1}}\prod_{j=0}^n(m+r+j-2+\sigma(j+1))\cdots(m+j+\sigma(j+1)).\\
\end{eqnarray*}

Therefore the multiset $\{i+\sigma(i)\mid i=1,\ldots,n+1\}$ is
determined by the expression for 
geometric genus, regarded as a polynomial in $m$.
This proves the claim and the theorem.
\end{pf}

\medskip

\noindent{\it Proof of the lemma.} Consider the quotient morphism 
$\pi :X^n\to X^{(n)}$. Then 
$H^i(X^{(n)},\bbC )=H^i(X^n,\bbC )^{S_n}$, where the action of $S_n$ on 
$H^i(X^n,\bbC )$ is twisted by the sign if $i$ is odd. Clearly, the 
$S_n$-action preserves the subspaces $H^{p,q}(X,\bbC )$. The embedding 
$H^*(X^{(n)},\bbC)\hookrightarrow H^*(X^n,\bbC)$ is a morphism of 
mixed Hodge structures; thus the Hodge structure on $H^*(X^{(n)},\bbC)$ is, 
in fact, pure and 
$$H^{p,q}(X^{(n)},\bbC)=H^{p,q}(X^n,\bbC)^{S_n}.$$
In particular,
$$H^{2n,0}(X^{(n)},\bbC)= Sym^nH^{2,0}(X,\bbC),$$
and
$$h^{2n,0}(X^{(n)})= \binom{r+n-1}{r-1}.$$
                               
It remains to prove that $P_g(X^{(n)})=h^{2n,0}(X^{(n)})$, i.e. 
$h^{2n,0}(\tilde{X})=h^{2n,0}(X^{(n)})$ for a resolution of singularities 
$\tilde{X}\to X$. 

Denote by $X^{[n]}$ the $n$th Hilbert scheme parametrizing closed 
zero-dimensional subschemes of length $n$ of $X$. It is known that 
$X^{[n]}$ is a smooth projective variety of dimension $2n$  and 
there exists a natural map $\pi: X^{[n]}\to X^{(n)}$ which is a resolution 
of singularities ([Na]). Let us show that $h^{2n,0}(X^{[n]})=
h^{2n,0}(X^{(n)})$. This follows immediately from the result of L.~Gottsche 
and W.~Soergel which computes the Hodge structure on the cohomology 
$H^\bullet(X^{[n]},\bbC)$. Let $P(n)$ be the set of partitions of $n$. 
We can write $\alpha \in P(n)$ as $n=\alpha _1 \cdot 1+\ldots +
\alpha _r\cdot r$. Put $|\alpha |=\sum \alpha _i$, 
$X^{(\alpha)}=X^{(\alpha _1)}\times \ldots \times X^{(\alpha _r)}$. 

\begin{thm}([GS]).  There exists a canonical isomorphism of mixed Hodge 
structures
$$H^{i+2n}(X^{[n]},\bbQ)\otimes \bbQ (n)=
\bigoplus_{\alpha \in P(n)}H^{i+2|\alpha|}(X^{(\alpha )},\bbQ)
\otimes \bbQ (|\alpha|).$$
\end{thm}

Here $\bbQ(1)$ is the one-dimensional Hodge structure of weight $-2$ and bidegree $(-1,-1)$.
The proof of this theorem is based on the decomposition theorem for 
mixed Hodge modules. 

Note that $\dim X^{(\alpha )}=2|\alpha |$. Since the Hodge structure 
$H^\bullet (X^{(\alpha )},\bbC)$ is a substructure of 
$H^\bullet (X^{\alpha _1}\times \ldots \times X^{\alpha _r},\bbC)$, it does 
not contain a summand $H^{p,q}$, unless $|p-q|\leq 2|\alpha|$. The same 
remains true after a twist by a tensor power of $\bbQ (1)$. Therefore the only summand on 
the RHS which contributes to $H^{2n,0}(X^{[n]},\bbQ)\otimes\bbQ(n)$ is 
$H^{0+2n}(X^{(n)},\bbQ)\otimes \bbQ (n)$. This proves the lemma.

\end{document}